\documentclass[12pt]{amsart}
\usepackage{fullpage}
\usepackage{fancyhdr}
\usepackage{latexsym,amsmath,amssymb,textcomp}
\usepackage{fancyhdr}
\usepackage{geometry}
\everymath{\displaystyle}
\newtheorem{theorem}{Theorem}[section]
\newtheorem{lemma}[theorem]{Lemma}
\newtheorem{proposition}[theorem]{Proposition}
\newtheorem{corollary}[theorem]{Corollary}
\newtheorem{definition}[theorem]{Definition\rm}
\newtheorem{remark}{\it \bf Remark\/}
\newtheorem{example}[theorem]{Example}

\newcommand{\mylabel}[1]{\label{#1}}

\newcommand{\aatop}[2]{\genfrac{}{}{0pt}{}{#1}{#2}}

\begin{document}

\title[Chern Invaraints]{THE CHERN INVARIANTS FOR PARABOLIC BUNDLES AT MULTIPLE POINTS}

\subjclass[2000]{Primary 14F05; Secondary 19L10, 14C17}

\keywords{Parabolic bundle, Chern character, Riemann-Roch, Gysin map, Bogomolov-Giesker Inequality, Birational transformation, Exceptional divisor, Elementary transformation}

\begin{abstract}
If $D\subset X$ is a curve with multiple points in a surface, a parabolic bundle
defined on $(X,D)$ away from the singularities can be extended in several ways
to a parabolic bundle on a resolution of singularities. We investigate the possible
parabolic Chern classes for these extensions.
\end{abstract}

\maketitle

\section{Introduction}
\mylabel{sec-intro}

Suppose $\check{X}$ is a smooth surface and
$\check{D}=\check{D}_1+\ldots + \check{D}_k$ is a divisor with each
$\check{D}_i$ smooth. Suppose $\check{E}$ is a bundle
provided with filtrations $\check{F}^i_{\cdot}$ along the $\check{D}_i$, and parabolic weights
$\alpha ^i_{\cdot}$. If $\check{D}$ has normal crossings, this defines a locally abelian
parabolic bundle on $(\check{X},\check{D})$ and the parabolic
Chern classes have been calculated as explained in the previous part.

Suppose that the singularities of $\check{D}$ contain some points of higher multiplicity.
For the present work we assume that these are as easy as possible, namely several
smooth branches passing through a single point with distinct tangent directions.
The first basic case is a triple point.

Let $\varphi : X \rightarrow \check{X}$ denote this
birational transformation, and let $D_i\subset X$ denote the strict transforms
of the $\check{D}_i$. Assuming for simplicity that there is a single multiple
point, denote by $D_0$ the exceptional divisor. Now $D=D_0+\cdots + D_k$ is
a divisor with normal crossings. Suppose $E$ is a vector bundle on $X$
with
$$
E|_{X-D_0} = \varphi ^{\ast}(\check{E})|_{X-D_0}.
$$
The filtrations $\check{F}^i_{\cdot}$
induce filtrations of $\varphi ^{\ast}(\check{E})|_{D_i}$ and hence of
$E|_{D_i-D_i\cap D_0}$, which then extend uniquely to filtrations $F^i_{\cdot}$ of
$E|_{D_i}$. Associate to these filtrations the same parabolic weights as before.

Up until now we have already made a choice of extension of the bundle
$E$. Choose furthermore a filtration $F^0_{\cdot}$ of $E|_{D_0}$
and parabolic weights associated to $D_0$. Having made these choices we get a parabolic
bundle on the normal crossings divisor $(X,D)$, which determines parabolic Chern classes.
We are particularly interested in the invariant $\Delta$ which combines $c_1$ and $c_2$
in such a way as to be invariant by tensoring with a line bundle.

The goal of this paper is to provide a convenient calculation of $\Delta$ and then
investigate its dependence on the choices which have been made above. In particular we would
like to show that $\Delta$ achieves its minimum and calculate this minimum, which
can be thought of as the Chern invariant associated to the original parabolic structure
on the multiple point singularity $(\check{X},\check{D})$.

The main difficulty is to understand the possible choices for $E$. For this
we use the technical of Ballico-Gasparim \cite{Ballico} \cite{BallicoGasparim1} \cite{BallicoGasparim2}.

\section{Calculating the invariant $\Delta$ of a locally abelian parabolic bundle}

Recall from \cite{Taher} formulas for the parabolic first, second Chern characters of a locally abelian parabolic bundle $E$
in codimension one and two, $ch_{1}^{Par}(E)$, and $ch_{2}^{Par}(E)$.

Let $X$ be a smooth projective variety over an algebraically closed
field of characteristic zero  and let $D$ be a strict normal crossings divisor
on $X$. Write $D = D_{1} + ... + D_{n}$ where $D_{i}$ are the irreducible
smooth components, meeting transversally. We sometimes denote by $\mathcal{S}:= \{ 1,\ldots , n\}$
the set of indices for components of the divisor $D$.

For $i = 1, ..., n$,
let $\Sigma_{i}$ be finite linearly ordered sets
with notations $ \eta_{i} \leq ... \leq \sigma \leq \sigma ' \leq \sigma ''
\leq ... \leq \tau_{i}$ where $\eta_{i}$ is  the smallest element of
$\Sigma_{i}$ and $\tau_{i}$ the greatest element of $\Sigma_{i}$.

Let $\Sigma '_{i}$ be the set of connections between the $\sigma$'s i.e
$$
\Sigma '_{i} = \{(\sigma, \sigma '),\ s.t\  \sigma < \sigma '\  and\
there \ exist \  no \ \sigma ''\  with\  \sigma < \sigma ''< \sigma'
\}.
$$
Consider the {\em tread functions} $ m_{+} :\Sigma'_i\rightarrow
\Sigma_i$ and $m_{-} : \Sigma'_i \rightarrow \Sigma_i$  if $\lambda =
(\sigma, \sigma') \in \Sigma'_{i} $ then $ \sigma = m_{-}(\lambda), \sigma' =
m_{+}(\lambda)$.
In the other direction, consider the {\em riser functions} $C_{+} : \Sigma _i- \{\tau _i\}
\rightarrow \Sigma'_i $ and $C_{-} : \Sigma_i - \{\eta_i\} \rightarrow
\Sigma '_i $ such that $C_{+}(\sigma) = (\sigma, \sigma ')$ where $
\sigma ' > \sigma$ the next element and $C_{-}(\sigma) = (\sigma '',
\sigma)$ where $\sigma '' < \sigma$ the next smaller element.
$$
$$
For any parabolic bundle $E$  in codimension one, and two,  the parabolic first, second  Chern characters
$ch_{1}^{Par}(E), and \  ch_{2}^{Par}(E),\ $  are obtained as follows:
$$
$$
$\centerdot \  ch^{Par}_{1}(E) := \  ch_{1}^{Vb}(E) \ - \ \sum_{i_{1} \in \mathcal{S}}\sum_{\lambda_{i_{1}} \in
\Sigma '_{i_{1}}}\alpha_{i_{1}}(\lambda_{i_{1}}).rank(Gr^{i_{1}}_{\lambda_{i_{1}}}).[D_{i_{1}}]$
$$
$$
$\centerdot \ ch^{Par}_{2}(E) := \  ch^{Vb}_{2}(E) \ -\  \sum_{i_{1} \in \mathcal{S}}\sum_{\lambda_{i_{1}} \in \sum'_{i_{1}}}
\alpha_{i_{1}}(\lambda_{i_{1}}).(\xi_{i_{1}})_{\star}\left(c_{1}^{D_{i_{1}}}(Gr^{i_{1}}_{\lambda_{i_{1}}})\right) $
$$
$$
$\hspace{2.2cm} +  \ \frac{1}{2} \ \sum_{i_{1} \in \mathcal{S}}\sum_{\lambda_{i_{1}} \in \sum'_{i_{1}}}
\alpha^{2}_{i_{1}}(\lambda_{i_{1}}).rank(Gr^{i_{1}}_{\lambda_{i_{1}}}).[D_{i_{1}}]^{2}$
$$
$$
$ \hspace{2.2cm} +  \ \frac{1}{2}\sum_{i_{1} \neq i_{2}}
\sum_{\aatop{\lambda_{i_{1}} }{ \lambda_{i_{2}}}}  \sum_{p \in Irr(D_{i_{1}} \cap D_{i_{2}})}
\alpha_{i_{1}}(\lambda_{i_{1}}).\alpha_{i_{2}}(\lambda_{i_{2}}).rank_{p}(Gr^{i_{1}, i_{2}}_{\lambda_{i_{1}}, \lambda_{i_{2}}}).[D_{p}].$
$$
$$
Where:
$\\$

$\centerdot \  ch_{1}^{Vb}(E), ch^{Vb}_{2}(E)$ denotes the first, second, Chern character of vector bundles $E$.

$\centerdot \  Irr(D_{I})$ denotes the set of the irreducible components of $D_I:=D_{i_{1}} \cap D_{i_{2}}  \cap ... \cap D_{i_{q}}$.

$\centerdot \  \xi_{I}$ denotes the closed immersion $D_{I} \longrightarrow X$, and $\xi_{I,{\star}} : A^k(D_{I}) \longrightarrow A^{k+q}(X)$ denotes  the
associated Gysin map.

$\centerdot$ \   Let $p$ be an element of $Irr(D_{i} \cap D_{j}).$  Then $rank_{p}(Gr^{I}_{\lambda})$ denotes the rank of $Gr^{I}_{\lambda}$ as
an $\mathcal{O}_{p}$-module.

$\centerdot  \  [D_{i_j}] \in A^1(X)\otimes \mathbb{Q}$,  and $[D_{p}] \in A^2(X)\otimes \mathbb{Q}$ denote the cycle classes given by $D_{i_j}$
and $D_{p}$ respectively.
\begin{definition}
Let $Gr^{i_{1}}_{\lambda_{i_{1}}}$ be a bundle over $D_{i_1}$. Define the degree of $Gr^{i_{1}}_{\lambda_{i_{1}}}$ to be
$$
\emph{deg} (Gr^{i_{1}}_{\lambda_{i_{1}}}) := (\xi_{i_{1}})_{\star}\left(c_{1}^{D_{i_{1}}}(Gr^{i_{1}}_{\lambda_{i_{1}}})\right).
$$
\end{definition}


\begin{definition}
The invariant $\Delta$, which is a normalized version of $c_2$ designed to be independent of tensorization by
line bundles. It is defined by
$$
\Delta = c_{2} - \frac{r - 1}{2r}c_{1}^{2}.
$$
Recall that:   $ch_{2} = \frac{1}{2}c_{1}^{2} - c_{2}$
$\Longrightarrow$  $c_{2} = \frac{1}{2}c_{1}^{2} - ch_{2}$.  Therefore
$$
\Delta = \frac{1}{2}c_{1}^{2} -  ch_{2} -  \frac{1}{2}c_{1}^{2}  +  \frac{1}{2}c_{1}^{2}  +  \frac{1}{2r}c_{1}^{2}  =  \frac{1}{2r}ch_{1}^{2} - ch_{2}.
$$
\end{definition}
Then $\Delta^{Par}(E) = \frac{1}{2r}ch^{Par}_{1}(E)^{2} - ch^{Par}_{2}(E) $
$$
$$
$ = \frac{1}{2r} \left[ch^{Vb}_{1}(E) \ - \ \sum_{i_{1} \in \mathcal{S}}\sum_{\lambda_{i_{1}} \in
\Sigma '_{i_{1}}}\alpha_{i_{1}}(\lambda_{i_{1}}).rank(Gr^{i_{1}}_{\lambda_{i_{1}}}).[D_{i_{1}}]\right]^{2}$
$$
$$
$ -\  ch^{Vb}_{2}(E)  \ +\  \sum_{i_{1} \in \mathcal{S}}\sum_{\lambda_{i_{1}} \in \sum'_{i_{1}}} \alpha_{i_{1}}(\lambda_{i_{1}}).\emph{deg}
(Gr^{i_{1}}_{\lambda_{i_{1}}}) $
$$
$$
$\hspace{2.2cm} - \ \frac{1}{2} \ \sum_{i_{1} \in \mathcal{S}}\sum_{\lambda_{i_{1}} \in \sum'_{i_{1}}}
\alpha^{2}_{i_{1}}(\lambda_{i_{1}}).rank(Gr^{i_{1}}_{\lambda_{i_{1}}}).[D_{i_{1}}]^{2}$
$$
$$
$ \hspace{2.2cm} - \ \frac{1}{2}\sum_{i_{1} \neq i_{2}}
\sum_{\aatop{\lambda_{i_{1}} }{ \lambda_{i_{2}}}}  \sum_{p \in Irr(D_{i_{1}} \cap D_{i_{2}})}
\alpha_{i_{1}}(\lambda_{i_{1}}).\alpha_{i_{2}}(\lambda_{i_{2}}).rank_{p}(Gr^{i_{1}, i_{2}}_{\lambda_{i_{1}}, \lambda_{i_{2}}}).[D_{p}].$
$$
$$
$\hspace{2.2cm} = \frac{1}{2r}ch^{Vb}_{1}(E)^{2} $
$$
$$
$\hspace{2.2cm} - \ \frac{1}{r} \left[ch^{Vb}_{1}(E) \right]
\cdot
\sum_{i_{1} \in \mathcal{S}}\sum_{\lambda_{i_{1}} \in
\Sigma '_{i_{1}}}\alpha_{i_{1}}(\lambda_{i_{1}}).rank(Gr^{i_{1}}_{\lambda_{i_{1}}}).[D_{i_{1}}]$
$$
$$
$ \hspace{2.2cm} + \ \frac{1}{2r}\sum_{i_{1} \neq i_{2}}
\sum_{\aatop{\lambda_{i_{1}} }{ \lambda_{i_{2}}}}  \sum_{p \in Irr(D_{i_{1}} \cap D_{i_{2}})} \alpha_{i_{1}}(\lambda_{i_{1}}).\alpha_{i_{2}}(\lambda_{i_{2}}).
rank(Gr^{i_{1}}_{\lambda_{i_{1}}})rank(Gr^{i_{2}}_{\lambda_{i_{2}}}).[D_{p}].$
$$
$$
$\hspace{2.2cm} + \ \frac{1}{2r}\sum_{i_{1} \in \mathcal{S}}\sum_{\lambda_{i_{1}} \in
\Sigma '_{i_{1}}}\sum_{\lambda'_{i_{1}} \in
\Sigma '_{i_{1}}}\alpha_{i_{1}}(\lambda_{i_{1}}).\alpha_{i_{1}}(\lambda'_{i_{1}}).
rank(Gr^{i_{1}}_{\lambda_{i_{1}}}).rank(Gr^{i_{1}}_{\lambda'_{i_{1}}}).[D_{i_{1}}]^{2}$
$$
$$
$\hspace{2.2cm} - \  ch^{Vb}_{2}(E)  \ +\  \sum_{i_{1} \in \mathcal{S}}\sum_{\lambda_{i_{1}} \in \sum'_{i_{1}}} \alpha_{i_{1}}(\lambda_{i_{1}}).\emph{deg}
(Gr^{i_{1}}_{\lambda_{i_{1}}}) $
$$
$$
$\hspace{2.2cm} - \ \frac{1}{2} \ \sum_{i_{1} \in \mathcal{S}}\sum_{\lambda_{i_{1}} \in \sum'_{i_{1}}}
\alpha^{2}_{i_{1}}(\lambda_{i_{1}}).rank(Gr^{i_{1}}_{\lambda_{i_{1}}}).[D_{i_{1}}]^{2}$
$$
$$
$ \hspace{2.2cm} - \ \frac{1}{2}\sum_{i_{1} \neq i_{2}}
\sum_{\aatop{\lambda_{i_{1}} }{ \lambda_{i_{2}}}}  \sum_{p \in Irr(D_{i_{1}} \cap D_{i_{2}})}
\alpha_{i_{1}}(\lambda_{i_{1}}).\alpha_{i_{2}}(\lambda_{i_{2}}).rank_{p}(Gr^{i_{1}, i_{2}}_{\lambda_{i_{1}}, \lambda_{i_{2}}}).[D_{p}].$

\begin{proposition}
 $\Delta^{Par}(E) = \Delta^{Vb}(E) $
 $$
 $$
 $\hspace{2.2cm} - \ \frac{1}{r} ch^{Vb}_{1}(E)
\cdot
\sum_{i_{1} \in \mathcal{S}}\sum_{\lambda_{i_{1}} \in
\Sigma '_{i_{1}}}\alpha_{i_{1}}(\lambda_{i_{1}}).rank(Gr^{i_{1}}_{\lambda_{i_{1}}}).[D_{i_{1}}]$
$$
$$
$ \hspace{2.2cm} + \ \frac{1}{2r}\sum_{i_{1} \neq i_{2}}
\sum_{\aatop{\lambda_{i_{1}} }{ \lambda_{i_{2}}}}  \sum_{p \in Irr(D_{i_{1}} \cap D_{i_{2}})} \alpha_{i_{1}}(\lambda_{i_{1}}).\alpha_{i_{2}}(\lambda_{i_{2}}).
rank(Gr^{i_{1}}_{\lambda_{i_{1}}})rank(Gr^{i_{2}}_{\lambda_{i_{2}}})[D_{p}].$
$$
$$
$\hspace{2.2cm} + \ \frac{1}{2r}\sum_{i_{1} \in \mathcal{S}}\sum_{\lambda_{i_{1}} \in
\Sigma '_{i_{1}}}\sum_{\lambda'_{i_{1}} \in
\Sigma '_{i_{1}}}\alpha_{i_{1}}(\lambda_{i_{1}}).\alpha_{i_{1}}(\lambda'_{i_{1}}).
rank(Gr^{i_{1}}_{\lambda_{i_{1}}}).rank(Gr^{i_{1}}_{\lambda'_{i_{1}}}).[D_{i_{1}}]^{2}$
$$
$$
$\hspace{2.2cm} + \  \sum_{i_{1} \in \mathcal{S}}\sum_{\lambda_{i_{1}} \in \sum'_{i_{1}}} \alpha_{i_{1}}(\lambda_{i_{1}}).\emph{deg}
(Gr^{i_{1}}_{\lambda_{i_{1}}}) $
$$
$$
$\hspace{2.2cm} - \ \frac{1}{2} \ \sum_{i_{1} \in \mathcal{S}}\sum_{\lambda_{i_{1}} \in \sum'_{i_{1}}}
\alpha^{2}_{i_{1}}(\lambda_{i_{1}}).rank(Gr^{i_{1}}_{\lambda_{i_{1}}}).[D_{i_{1}}]^{2}$
$$
$$
$ \hspace{2.2cm} - \ \frac{1}{2}\sum_{i_{1} \neq i_{2}}
\sum_{\aatop{\lambda_{i_{1}} }{ \lambda_{i_{2}}}}  \sum_{y \in Irr(D_{i_{1}} \cap D_{i_{2}})}
\alpha_{i_{1}}(\lambda_{i_{1}}).\alpha_{i_{2}}(\lambda_{i_{2}}).rank_{y}(Gr^{i_{1}, i_{2}}_{\lambda_{i_{1}}, \lambda_{i_{2}}}).[D_{y}].$
\end{proposition}
\section{Parabolic bundles with full flags}

\mylabel{sec-parasurface}
We use the fact that $X$ is a surface to simplify the above expressions, by assuming that the parabolic filtrations are
full flags.

\begin{proposition}
If $E'$ is a locally free sheaf over $X - \lbrace P \rbrace$,  then $\exists !$ extension to a locally free sheaf $E$ over $X$ s.t
$E'\vert_{X - \lbrace P\rbrace} = E$.
\end{proposition}
\begin{proposition}
If we have a strict sub-bundle of $E\vert D_{i} - \lbrace P \rbrace$ then $\exists !$ extension to a strict sub-bundle of $E\vert D_{i}$.
\end{proposition}
\begin{remark}
It follows from
these propositions that if $(E', F'^{i}_{\alpha_{i}})$ is a parabolic structure over $(X - \lbrace P \rbrace, D - \lbrace P \rbrace)$, then
we obtain a bundle $E$ over $X$ with the filtrations $\lbrace F'^{i}_{\alpha_{i}} \rbrace$ of $E\vert_{D_{i}}$ by a strict sub-bundles.
\end{remark}
\begin{definition}
$Hom(\mathcal{O}_{\mathbb{P}^{1}}(-m_{i}), \mathcal{O}_{\mathbb{P}^{1}}^{\oplus2}) = H^{0}(\mathbb{P}^{1}, \mathcal{O}(-m_{i})\otimes \mathcal{O})^{\oplus2} =
H^{0}(\mathbb{P}^{1}, \mathcal{O}(m_{i}))^{\oplus2}.$
\end{definition}
For example, subbundles of a rank two trivial bundle may be expressed very explicitly.

\begin{proposition}
Consider the two polynomials $(A_{i}, B_{i}) \in H^{0}(\mathbb{P}^{1}, \mathcal{O}(m_{i}))^{\oplus2}$,
the sub-sheaves are saturated iff $m_{i} = ((max(deg(A_{i}), deg(B_{i}))$ and $(A_{i}, B_{i}) = 1$. Then there is an  isomorphism
$$
(A_{i}, B_{i}) : \mathcal{O}_{\mathbb{P}^{1}}(-m_{i}) \longrightarrow \mathcal{O}_{\mathbb{P}^{1}}^{\oplus2}
$$
\end{proposition}
\begin{lemma}
$\forall$   $0 \subseteq F^{i}_{\sigma_{i}}  \subseteq F^{i}_{\sigma'_{i}} \subseteq  .... \subseteq F^{i}_{\tau_{i}} \subseteq E\vert_{D_{i}}$,  $\exists$
complet flags  $0 \subseteq \widehat{F}_{1} \subseteq ... \subseteq \widehat{F}_{r}= E\vert_{D_{i}}$ s.t $F^{i}_{\sigma_{i}} = \widehat{F}_{k(\sigma_{i})}$.
where $\forall$ $\sigma_{i} \in \Sigma_{i}$
we have  $k(i) \in \lbrace 0, 1, ..., r \rbrace$ then $\exists$ $k : \Sigma_{i} \longrightarrow \lbrace 0, 1, ..., r\rbrace$
s.t $k = rank(F^{i}_{\sigma_{i}})$.
\end{lemma}

In view of this lemma, we will now suppose that all the filtrations are complete flags.
The weights should then form an increasing sequence but not necessarily strictly
increasing.

In particular we will change notation
and denote the filtration of $E\vert_{D_{i}}$ by
$$
0= F^i_0 \subseteq F^i_1\subseteq ....\subseteq F^i_r = E\vert_{D_{i}} .
$$
In this case $\Sigma _i = \{ \sigma _i^0,\ldots , \sigma _i^r\}$ and
$\Sigma '_i = \{ \lambda _i^1,\ldots , \lambda _i^r\}$.
These sets now have the same number of elements for each $i$ so we can return to a numerical indexation.
We denote
$$
Gr^i_k(E\vert_{D_{i}}) := Gr^i_{\lambda _i^k}(E\vert_{D_{i}}) = F^i_k/F^i_{k-1}.
$$
\begin{proposition}
$rank(Gr^{i}_{\lambda_{i}^{}}) = 1$ $ \Longleftrightarrow$ the filtrations
$ F^i_{\sigma ^1_i} \leq
F^i_{\sigma ^2_i} \leq ...\leq F^i_{\sigma ^r_i}$  are complet flags  for $i = 1,2,...,n$.
\end{proposition}
Since we are on a surface,
$D_i\cap D_j$ is a finite collection of points. At each point $P\in D_i\cap D_j$ we have two filtrations of $E_P$
coming from the parabolic filtrations along $D_i$ and $D_j$. We are now assuming that they are both complete flags.
The incidence relationship between these filtrations is therefore encoded by a permutation.

\begin{lemma} $\forall$ $k$ $\exists !$ $k' \in  \lbrace 1, ..., r \rbrace$  s.t $\hspace{0.1cm}$ $rank \left(Gr^i_{k} Gr^i_{k'}(E_{P})\right)= \dfrac{F_{k}^{i}
\cap F_{k'}^{j}}{F_{k-1}^{j} \cap F_{k'}^{j} + F_{k}^{i} \cap F_{k'-1}^{i}} = 1$.
\end{lemma}
\begin{definition}
$\forall$  $P \in D_{i} \cap D_{j}$ define the permutation
$\sigma( P, i, j) : \lbrace 1, ..., r \rbrace= \Sigma'_{i} \longrightarrow \lbrace 1, ..., r \rbrace $ which sends  $k \in \lbrace 1, ..., r \rbrace $ to $
\sigma( P, i, j)(k) = k'$ where $k'$ is the unique index given in the previous lemma.
\end{definition}
\begin{lemma}
$\forall$ $k$ if $k'' \neq \sigma(P, i, j)(k)$ then $ \hspace{0.1cm}$ $rank \left(Gr^{i}_{k} Gr^{i}_{k''}(E_{P})\right) = 0$.
\end{lemma}

Since the filtrations are full flags, there are $r$ different indices $\lambda _i^1,\ldots , \lambda _i^r$ for each divisor $D_i$.
We introduce the notation $\alpha (D_i,k):=\alpha _i(\lambda _i^k)$.

With this notation we obtain the following expression for the term involving $Gr^{i_{1}, i_{2}}_{\lambda_{i_{1}}, \lambda_{i_{2}}}$:
$$
$$
$ \hspace{2.2cm} - \ \frac{1}{2}\sum_{i_{1} \neq i_{2}}
\sum_{\aatop{\lambda_{i_{1}} }{ \lambda_{i_{2}}}}  \sum_{p \in Irr(D_{i_{1}} \cap D_{i_{2}})}
\alpha_{i_{1}}(\lambda_{i_{1}}).\alpha_{i_{2}}(\lambda_{i_{2}}).rank_{y}(Gr^{i_{1}, i_{2}}_{\lambda_{i_{1}}, \lambda_{i_{2}}}).[y].$
$$
$$
$= \hspace{2.2cm} - \ \frac{1}{2}\sum_{i \neq j}
\sum_{k=1}^r  \sum_{y \in Irr(D_{i} \cap D_{j})} \alpha (D_i,k).\alpha (D_j,\sigma (y,i,j)(k)).[y].$
$$
$$
On the other hand, all ranks of the graded pieces $Gr(D_i,k):=Gr^i_{\lambda ^k_i}$ are equal to $1$.
They are line bundles on $D_i$.
$\\$
\begin{definition}
Suppose $i \neq 0$. Let $Gr(D_i, k)$ are line bundles over $D_i$. Then we define the $deg Gr(D_i, k)$ to be:
$$
\emph{deg}\left(Gr(D_i, k)\right) = (\xi_{i})_{\star}\left(c_{1}^{D_{i}}(Gr(D_i,k))\right).
$$
\end{definition}

We can now rewrite the statement of Proposition.

\begin{proposition}
$\Delta^{Par}(E) = \Delta^{Vb}(E) $
$$
$$
$\hspace{2.2cm} - \ \frac{1}{r} ch^{Vb}_{1}(E)
\cdot
\sum_{i \in \mathcal{S}}\sum_{k=1}^r\alpha (D_i,k).[D_i]$
$$
$$
$ \hspace{2.2cm} + \ \frac{1}{2r}\sum_{i\neq j}
\sum_{k,l\in [1,r]}  \sum_{y \in Irr(D_{i} \cap D_{j})} \alpha (D_i,k).\alpha (D_j,l)[y].$
$$
$$
$\hspace{2.2cm} + \ \frac{1}{2r}\sum_{i \in \mathcal{S}}\sum_{k,l\in [1,r]}
\alpha (D_i,k).\alpha (D_i,l)
.[D_{i}]^{2}$
$$
$$
$\hspace{2.2cm} + \  \sum_{i \in \mathcal{S}}\sum_{k=1}^r \alpha (D_i,k).\emph{deg}\left(Gr(D_i, k)\right) $
$$
$$
$\hspace{2.2cm} - \ \frac{1}{2} \ \sum_{i\in \mathcal{S}}\sum_{k=1}^r
\alpha (D_i,k)^{2}.[D_{i}]^{2}$
$$
$$
$\hspace{2.2cm} - \ \frac{1}{2}\sum_{i \neq j}
\sum_{k=1}^r  \sum_{y \in Irr(D_{i} \cap D_{j})} \alpha (D_i,k).\alpha (D_j,\sigma (y,i,j)(k)).[y].$
\end{proposition}

We have $\alpha(D_i, k) \in[-1, 0]$, define $\alpha^{tot}(D_i) := \sum_{k=1}^{r} \alpha(D_i, k)$.
With this notation,

$\Delta^{Par}(E) = \Delta^{Vb}(E) $
$$
$$
$\hspace{2.2cm} - \ \frac{1}{r} ch^{Vb}_{1}(E)
\cdot
\sum_{i \in \mathcal{S}}\alpha^{tot}(D_i)[D_i]$
$$
$$
$ \hspace{2.2cm} + \ \frac{1}{2r}\sum_{i\neq j}
\alpha^{tot}(D_i)\alpha^{tot}(D_j)[D_i\cap D_j]$
$$
$$
$\hspace{2.2cm} + \ \frac{1}{2r}\sum_{i \in \mathcal{S}}
\alpha^{tot}(D_i)^2
.[D_{i}]^{2}$
$$
$$
$\hspace{2.2cm} + \  \sum_{i \in \mathcal{S}}\sum_{k=1}^r \alpha (D_i,k). $
$$
$$
$\hspace{2.2cm} - \ \frac{1}{2} \ \sum_{i\in \mathcal{S}}\sum_{k=1}^r
\alpha (D_i,k)^{2}.[D_{i}]^{2}$
$$
$$
$\hspace{2.2cm} - \ \frac{1}{2}\sum_{i \neq j}
\sum_{k=1}^r  \sum_{y \in Irr(D_{i} \cap D_{j})} \alpha (D_i,k).\alpha (D_j,\sigma (y,i,j)(k)).[y].$
$$
$$
For more simplification of $\Delta^{Par}(E)$, define $\beta$ such that
$$
$$
$\beta(D_i, k) := \alpha(D_i, k) - \frac{\alpha^{tot}(D_i)}{r} \hspace{0.2cm} \Longrightarrow \hspace{0.2cm} \alpha(D_i, k) = \beta(D_i, k) +
\frac{\alpha^{tot}(D_i)}{r}$.
$$
$$
\begin{remark}
We remark that $\sum_{k=1}^{r}\beta(D_i, k) = 0$. Hence
$$
\sum_{k=1}^{r}\alpha (D_i, k)^2 = \sum_{k=1}^{r}\beta(D_i, k)^2 +
\frac{\alpha ^{tot}(D_i)^2}{r},
$$
and for $i\neq j$ and $y \in Irr(D_{i} \cap D_{j})$,
$$
\sum_{k=1}^r  \alpha (D_i,k).\alpha (D_j,\sigma (y,i,j)(k)) =
\sum_{k=1}^r  \alpha (D_i,k).\alpha (D_j,\sigma (y,i,j)(k)) +
\frac{ \alpha ^{tot}(D_i).\alpha ^{tot}(D_j)}{r}.
$$
Furthermore note that
$$
\sum_{k=1}^r (\xi_{i})_{\star}\left( c_{1}^{D_{i}}(Gr(D_i,k))\right) = (\xi_{i})_{\star}\left( c_{1}^{D_{i}}(E|_{D_i})\right) = c_1^{Vb}(E).[D_i]
$$
so
$$
\sum_{k=1}^r \alpha (D_i,k).\emph{deg}\left(Gr(D_i, k)\right)
=
\sum_{k=1}^r \beta (D_i,k).\emph{deg}\left(Gr(D_i, k)\right) +
\frac{\alpha ^{tot} (D_i) c_1^{Vb}(E).[D_i]}{r} .
$$
\end{remark}
Using these remarks and the previous formula we get
$\\$

$\Delta^{Par}(E) = \Delta^{Vb}(E) $
$$
$$
$\hspace{2.2cm} - \ \frac{1}{r} ch^{Vb}_{1}(E)
\cdot
\sum_{i \in \mathcal{S}}\alpha^{tot}(D_i)[D_i]$
$$
$$
$ \hspace{2.2cm} + \ \frac{1}{2r}\sum_{i\neq j}
\alpha^{tot}(D_i)\alpha^{tot}(D_j)[D_i\cap D_j]$
$$
$$
$\hspace{2.2cm} + \ \frac{1}{2r}\sum_{i \in \mathcal{S}}
\alpha^{tot}(D_i)^2
.[D_{i}]^{2}$
$$
$$
$\hspace{2.2cm} + \  \sum_{i \in \mathcal{S}}\sum_{k=1}^r \beta (D_i,k). \emph{deg}\left(Gr(D_i, k)\right)$
$$
$$
$\hspace{2.2cm} + \ \frac{1}{r}  \sum_{i \in \mathcal{S}}\alpha ^{tot} (D_i).c_1^{Vb}(E).[D_i] $
$$
$$
$\hspace{2.2cm} - \ \frac{1}{2} \ \sum_{i\in \mathcal{S}}\sum_{k=1}^r
\beta (D_i,k)^{2}.[D_{i}]^{2}$
$$
$$
$\hspace{2.2cm} - \ \frac{1}{2r} \ \sum_{i\in \mathcal{S}}
\alpha ^{tot}(D_i)^{2}.[D_{i}]^{2}$
$$
$$
$\hspace{2.2cm} - \ \frac{1}{2}\sum_{i \neq j}
\sum_{k=1}^r  \sum_{y \in Irr(D_{i} \cap D_{j})} \beta (D_i,k).\beta (D_j,\sigma (y,i,j)(k)).[y].$
$$
$$
$\hspace{2.2cm} - \ \frac{1}{2r}\sum_{i \neq j}
\sum_{k=1}^r  \sum_{y \in Irr(D_{i} \cap D_{j})} \alpha ^{tot}(D_i).\alpha ^{tot}(D_j).[y].$
$$
$$
The terms containing $\alpha ^{tot}(D_i)$ all cancel out, giving the following
formula.

\begin{proposition}
\mylabel{propepar}
$\Delta^{Par}(E) = \Delta^{Vb}(E) $
$$
$$
$\hspace{2.2cm} + \  \sum_{i \in \mathcal{S}}\sum_{k=1}^r \beta (D_i,k).\emph{deg}\left(Gr(D_i, k)\right) $
$$
$$
$\hspace{2.2cm} - \ \frac{1}{2} \ \sum_{i\in \mathcal{S}}\sum_{k=1}^r
\beta (D_i,k)^{2}.[D_{i}]^{2}$
$$
$$
$\hspace{2.2cm} - \ \frac{1}{2}\sum_{i \neq j}
\sum_{k=1}^r  \sum_{y \in Irr(D_{i} \cap D_{j})} \beta (D_i,k).\beta (D_j,\sigma (y,i,j)(k)).[y].$
\end{proposition}

The fact that $\Delta ^{Par}(E)$ is independent of $\alpha ^{tot}$ is the parabolic
version of the invariance of $\Delta$ under tensoring with line bundles. Even though
this is the theoretical explanation, for the proof it was more
convenient to calculate explicitly the formula and notice that the
terms containing $\alpha ^{tot}$ cancel out, than to try to compute the tensor
product with a parabolic line bundle.
\section{Resolution of singular divisors}

Now we can consider a more general situation, where $\Check{X}$ is a smooth projective surface but $D=\bigcup _{i=1}^n D_i$ is a divisor
which may have singularities worse than normal crossings. Let $\Check{P} = \lbrace \Check{P}_{1}, ..., \Check{P}_{r} \rbrace$ be a set of points. Assume that the points $\Check{P}_{j}$ are
crossing points of $\Check{D}_{i}$, and that they are general multiple points, that is through a crossing point $P_{j}$ we have divisors $\Check{D}_{i_1}, ..., \Check{D}_{i_m}$ which are
pairwise transverse. Assume that $\Check{D}$ has normal crossings outside of the set of points $\Check{P}$.
We choose an embedded resolution given by a sequence of blowing-ups
$\varphi:X\rightarrow \check{X}$ in $r$ points $\Check{P}_{1}, ..., \Check{P}_{r}$  and $P$ be the  exceptional divisor on $X$, note that
$P$ is a sum of disjoint exceptional components $P_{i} = \varphi^{-1}(\check{P}_{i})$ over the points $\check{P}_{i}$ respectively. The pullback divisor may be written as $D=D_1+\cdots +D_a +P_1+\cdots + P_b$ where $D_i$ is the strict transform of a component $\Check{D}_i$ of
the original divisor, and $P_j$ are the exceptional divisors.
\begin{definition}
Let $E$ be a bundle over $X$, and consider the inclusion $i : U \hookrightarrow X$ where $U = X - \bigcup_{i=1}^{k}{P_{i}}$ be a smooth connected
quasi-projective surface. Hence $P_{i} = \mathbb{P}^{1}$  and let the blowing-up  $\varphi : X \longrightarrow \check{X} $. Define $\check{E}$ as a unique bundle
over $\check{X}$  such that
$$
\check{E}|_{U} \cong E|_{U},
$$
$$
\check{E} \hspace{0.15cm}  is \hspace{0.15cm} locally \hspace{0.15cm} free.
$$
\end{definition}

This construction allows us to localize the contributions of the Chern classes of $E$ along the exceptional divisors,
by comparison with $\varphi ^{\ast} (\check{E})$.
\begin{definition}
Let $E$ be a bundle over $X$. Consider the inclusions $\varphi^{\star}\check{E} \hookrightarrow i_{\star}(E|_{U})$, where $i_{\star}(E|_{U})$ is a quasi-coherent
sheaves over $X$, and  $E \hookrightarrow i_{\star}(E|_{U})$, where $i : U \hookrightarrow \check{X}$. Define $E''$ to be the intersection of subsheaves
$\varphi^{\star}\check{E}$ and $E$ of $i_{\star}(E|_{U})$.
\end{definition}
\begin{lemma}
$E''$ is a free locally coherent sheaf.
\end{lemma}
\begin{definition}
Consider the two exact sequences
$$
0 \longrightarrow E'' \longrightarrow \varphi^{\star}\check{E} \longrightarrow Q' \longrightarrow 0
$$
$$
0 \longrightarrow E'' \longrightarrow E \longrightarrow Q \longrightarrow 0
$$
Let $E/E'' = Q = \bigoplus_{i=1}^{k} Q_{i}$ and $\varphi^{\star}\check{E}/E'' = Q' = \bigoplus_{i=1}^{k} Q'_{i}$. Define the local contribution to be,
$$
ch^{Vb}(E,P)_{loc} := ch^{Vb}(Q) - ch^{Vb}(Q')
$$
\end{definition}
\begin{proposition}
If  $\varphi : \check{X} \longrightarrow X$. Let  $P_{i}$ the blowing-up of $\check{P_{i}}$ where $P_{i}$ is the exceptional divisor for $i = 1, 2, ...k$. Then
$$
ch^{Vb}(E) = ch^{Vb}(\varphi^{\star}(\check{E})) + \sum_{i=1}^{k}ch^{Vb}(E, P_{i})_{loc}
$$
\end{proposition}
We have $ch^{Vb}_{1}(E) \in A^{1}(X)$. Let $\varphi^{\star} : A^{1}(\check{X}) \longrightarrow A^{1}(X)$; where

$A^{1}(X) =  A^{1}(\check{X}) \oplus \bigoplus_{i=1}^{k} .\mathbb{Z}.[P_{i}].$  We have

$P_{i}.\varphi^{\star}(\check{D}) = 0$ if $\check{D} \in A^{1}(\check{X})$ and $P_{i}.P_{j} = 0$ if $i \neq j$. Then

$ch^{Vb}_{1}(E) = \varphi^{\star}ch^{Vb}_{1}(\check{E}) + \sum_{i = 1}^{k}a_{i}[P_{i}] = \varphi^{\star}ch^{Vb}_{1}(\check{E}) + \sum_{i = 1}^{k}ch^{Vb}_{1}(E,
P_{i})_{loc}.$

When we take the square, the cross-terms are zero, indeed $ch^{Vb}_{1}(E, P_{i})_{loc}$ is a multiple of the divisor class $[P_i]$
but $[P_i].[P_j] = 0 $ for $i\neq j$, and $[P_i].\varphi ^{\ast}[C]=0$ for any divisor $C$ on $\Check{X}$. Therefore,

$ch^{Vb}_{1}(E)^{2} = \varphi^{\star}ch^{Vb}_{1}(\check{E})^{2} + \sum_{i = 1}^{k}a_{i}^{2}[P_{i}]^{2} = \varphi^{\star}ch^{Vb}_{1}(\check{E})^{2} + \sum_{i =
1}^{k}ch^{Vb}_{1}(E, P_{i})^{2}_{loc}$.

\begin{lemma}
If $L $ is a line bundle over $X$, then
$$
\Delta(E \otimes L) = \Delta(E).
$$
\end{lemma}
\section{Local Bogomolov-Gieseker inequality}

The classical {\em Bogomolov-Gieseker inequality} states that if $X$ is projective and $E$ is a semistable vector bundle then
$\Delta (E)\geq 0$. We will see that a local version holds; the first observation is that the invariant $\Delta$ can be localized,
even though it involves a quadratic term in $ch_1$.

\begin{definition}
$$
\Delta^{Vb}(E, P_{i})_{loc} := \frac{1}{2r}ch_{1}^{Vb}(E, P_{i})^{2}_{loc}- ch_{2}^{Vb}(E, P_{i})_{loc}
$$
\end{definition}
\begin{lemma}
If $L = \varphi^{\star} \check{L}(\sum b_{i}P_{j})$ is a line bundle over $X$. Then
$$
\Delta^{Vb}(E \otimes L; P_{i})_{loc} = \Delta^{Vb}(E, P_{i})_{loc}
$$
\end{lemma}
\begin{proposition}
$$
\Delta^{Vb}(E) = \varphi^{\star}\Delta^{Vb}(\check{E}) + \sum_{i=1}^{k}\Delta^{Vb}(E, P_{i})_{loc}
$$
\end{proposition}

In order to get a bound, the technique is
to apply the Grothendieck decomposition to analyse more closely the structure of $E$ near the
exceptional divisors $P_i$, following Ballico \cite{Ballico} and Ballico-Gasparim \cite{BallicoGasparim1} \cite{BallicoGasparim2}
and others.
\begin{theorem}
Every vector bundle $E$  on $\mathbf{P}^{1}$ is of the form  $\mathcal{O}(m_{{1}})^{r_1} \oplus \cdots \oplus \mathcal{O}(m_{{a}})^{r_a} = \bigoplus_{j =1}^{a}
\mathcal{O}(m_{{j}})^{r_j}$,  $m_{{1}} < \ldots  <m_{{r}}$ where $m_{{j}} \in \mathbb{Z}$, and the
$r_j$ are positive integers with $r_1+\ldots + r_a=r$.  This called the Grothendick decomposition and it is unique.
\end{theorem}

Apply this decomposition to the restriction of the bundle $E$ to each exceptional divisor $P_i\cong \mathbf{P}^{1}$.
Thus
$$
E|_{P_{i}} = \mathcal{O}(m_{i,{1}})^{r_{i,1}}\oplus ... \oplus \mathcal{O}(m_{i,{a_i}})^{r_{i,a_i}}
= \bigoplus_{j =1}^{a_i} \mathcal{O}(m_{i,{j}})^{r_{i,j}}
$$
with $m_{i,{1}}< \ldots <  m_{i,{a_i}}$.

\begin{proposition}
Let $E$ be a bundle over $X$,  we have,
$$
m_{i,{j}} \hspace{0.1cm} = 0 \Longleftrightarrow \hspace{0,2cm} E \cong \varphi^{\star}\check{E} \hspace{0.2cm},
$$
if $E' = E(\sum_{i} k_i.P_{i})$ then $m'_{i,{j}} = m_{i,{j}} -k_i$, therefore
$$
m_{i,{j}} = k_i \Longleftrightarrow E \cong (\varphi^{\star}\check{E})(-\sum_{i} k_i.P_{i}).
$$
In this case we say that $E$ is {\em pure}, it is equivalent to saying that $a_i=1$.
\end{proposition}

See Ballico-Gasparim \cite{BallicoGasparim1}.
\begin{definition}
Let $E$ be a non trivial bundle, and $E|_{P} = \mathcal{O}(m_1)^{r_1} \oplus ... \oplus \mathcal{O}(m_a)^{r_a}$ be the restriction of the bundle $E$, for $m_1 <
m_2 < ... < m_a$. We define
$$
min(E|_{P}) := m_1 \hspace{0.2cm}, max(E|_{P}) := m_a, \hspace{0.2cm} and  \hspace{0.2cm}\varphi(E) = max(E|_{P}) - min(E|_{P}).
$$
\end{definition}
\begin{remark}
If $\mu(E) = m_a - m_1 = 0$. Then $E|_{P} = \mathcal{O}_{\mathbb{P}^1}(m_1)^r$ ; $E = E^{\vee}(-m.P)$
\end{remark}
\begin{lemma}
\mylabel{minmax}
If we have an exact sequence of bundles over $\mathbb{P}^1$
$$
0 \longrightarrow U \longrightarrow V \longrightarrow W \longrightarrow 0
$$
then
$$
min(V) \geq min(min(U), min(W)),
$$
$$
max(V) \leq max(max(U), max(W)).
$$
\begin{proof}
Define
$$
max(U) = max\{n; \hspace{0.05cm} s.t \hspace{0.15cm} \exists \hspace{0.15cm} \mathcal{O}_{\mathbb{P}^1}(n) \rightarrow U  \hspace{0.05cm} non \hspace{0.05cm}
trivial \} = max \{ n; \hspace{0.05cm} s.t \hspace{0.15cm}H^0(U(-n)) \neq 0\}
$$
$$
min(U) = min\{n; \hspace{0.05cm} s.t \hspace{0.15cm} \exists \hspace{0.15cm} U \rightarrow \mathcal{O}_{\mathbb{P}^n}(n)  \hspace{0.05cm} non \hspace{0.05cm}
trivial \} = max \{ n; \hspace{0.05cm} s.t \hspace{0.15cm}H^0(U^\ast(n)) \neq 0\}
$$
then
$$
max(U) \leq max(max(U), max(W))
$$
$$
min(V) \geq min(min(U), min(W))
$$
\end{proof}
\end{lemma}
Now we concentrate on one of the exceptional divisors $P_i$
and supress the index $i$ from the notation.

Now for $1 \leq t \leq r$, suppose that $E|_P$ is not pure, and consider the exact sequence
$$
\begin{array}{ccl}
0  &  &  \\
\uparrow & & \\
Q & := & \mathcal{O}(m_1)^{r_1}  \\
\uparrow & & \\
E|_{P_{i}} & := & \mathcal{O}(m_1)^{r_1} \oplus \mathcal{O}(m_2)^{r_2} \oplus \cdots \oplus \mathcal{O}(m_a)^{r_a}\\
\uparrow & & \\
K &:= & \mathcal{O}(m_{2})^{r_2} \oplus \cdots \oplus \mathcal{O}(m_{a})^{r_a}\\
\uparrow & & \\
0& &
\end{array}
$$
\begin{definition}
Suppose $X$ and $D$ are smooth with $D \stackrel{i_\ast}{\hookrightarrow} X$. Let  $E$ be a  free locally bundle over $X$.
Suppose we have an exact sequence
$$
0 \longrightarrow K \longrightarrow E|_D \longrightarrow Q \longrightarrow 0
$$
where $Q = \mathcal{O}(m_{1})^{r_{1}}$ is called constant stabilizer.
Define $E'$  to be the elementary transformation of $E$ by
$$
E':= Ker(E \rightarrow i_{\ast} Q).
$$
Then the sequence
$$
0 \longrightarrow E' \longrightarrow E \longrightarrow i_{\ast}Q \longrightarrow 0.
$$
is exact.
\end{definition}
\begin{lemma}
\mylabel{elemexact}
We have an exact sequence
$$
0 \longrightarrow Q(-D) \longrightarrow E'\mid_P \longrightarrow K \longrightarrow 0.
$$
Then
$$
E'(U) = \{ S \hspace{0.05cm} \in \hspace{0.05cm} E(U) \hspace{0.05cm} s.t \hspace{0.05cm} S|_{(D \cap U)} \in K(D \cap U)\}.
$$
\end{lemma}
\begin{lemma}
$\mu(E') \leq \mu(E) -1$ (if $\mu(E) \geq 1$).
\begin{proof}
We have
$$
\mathcal{O}_P(-P) = \mathcal{O}_P(i)
$$
apply the Lemma \ref{elemexact} we get
$$
0 \longrightarrow \mathcal{O}(m_1 + 1)^{r_1} \longrightarrow E'|_P  \longrightarrow \bigoplus_{i = 2}^{a} \mathcal{O}(m_i)^{r_i} \longrightarrow 0
$$
apply Lemma \ref{minmax}, take $min = max = m_1 + 1$ and $min = m_2$ then $max \geq m_1 + 1$

$\Longrightarrow$ $min(E'\mid_P) \leq m_a$. Therefore
$$
\mu(E') \leq \mu(E) - 1.
$$
\end{proof}
\end{lemma}
\begin{lemma}
 $$
 ch\left(\mathcal{O}_P(m_1)\right) =(P +
(m_{1} + \frac{1}{2})).
 $$
\end{lemma}
\begin{proof}
We have
$$
\mathcal{O}(m_{1}) = \mathcal{O}(-m_{1}P),
$$
consider the exact sequence
$$
0\longrightarrow \mathcal{O}_{X}(-(m_{1} + 1)P) \longrightarrow \mathcal{O}_{X}(-m_{1}P) \longrightarrow \mathcal{O}_{P}(-m_{1}P) \longrightarrow 0
$$
then

$ch(\mathcal{O}_{P}(-m_{1} P)) = e^{-m_{1}P} - e^{(-m_{1} + 1)P} = (1 - m_{1}P + \frac{m_{1}^{2}}{2}P^{2}) - (P - m_{1}P^{2} - \frac{P^{2}}{2})$

$ = (P - (m_{1} + \frac{1}{2})P^2)$, but $P^{2} = -1$, then $ch(\mathcal{O}_{P}(-m_{1})) =  (P + (m_{1} + \frac{1}{2}))$.
\end{proof}
\begin{proposition}
We have

$$
ch^{Vb}_{1}(E') = ch^{Vb}_{1}(E) - r_{1}P  \hspace{0.3cm} and \hspace{0.3cm}  ch^{Vb}_{2}(E') = ch^{Vb}_{2}(E) - (m_{1} + \frac{1}{2})r_{1},
$$
then
$$
\Delta^{Vb}(E', P) = \frac{1}{2r}(ch^{Vb}_{1}(E) - r_{1}P)^{2} - ch^{Vb}_{2}(E) + (m_{1} + \frac{1}{2})r_{1},
$$
therefore
$$
\Delta^{Vb}(E', P)= \Delta^{Vb}(E, P) - \frac{r_{1}}{r}ch^{Vb}_{1}(E).P - \frac{r_{1}^{2}}{2r} + (m_{1} + \frac{1}{2})r_{1}.
$$
\end{proposition}

We can now calculate using the previous lemma.
$$
$$
$E|_{P} = \mathcal{O}(m_{1})^{r_{1}} \oplus  \mathcal{O}(m_{2})^{r_{2}} \oplus ... \oplus  \mathcal{O}(m_{k})^{r_{k}}$,
where $ \sum_{i=1}^{k} r_{i} = r$, we have

$ch_{1}(E).P = \xi_{P,\star} \left(ch_{1}(E\mid_{P}) \right) = ch_{1} \left(\bigoplus_{i = 1}^{k}\mathcal{O}(m_{i}) \right) = \sum_{i=1}^{k} m_{i}r_{i}$, then

$\Delta^{Vb}(E', P) = \Delta^{Vb}(E, P) - \frac{r_{1}}{r}\sum_{i=1}^{k} m_{i}r_{i} - \frac{r_{1}^{2}}{2r} + m_{1}r_{1} + \frac{r_{1}}{2} = \Delta(E) -
\frac{1}{r} \mathcal{A}$, where

$\mathcal{A} = \sum_{i=2}^{k} m_{i} r_{i} r_{1} + m_{1}r_{1}^{2} + \frac{1}{2} r_{1}^{2} - (m_{1} + \frac{1}{2}) r_{1}r  \hspace{0.4cm}  for \hspace{0.4cm} r = r_{1} +
.... + r_{k}$

$ = \sum_{i=2}^{k}m_{i}  r_{i} r_{1} + m_{1}r_{1}^{2} + \frac{1}{2} r_{1}^{2} - (m_{1} + \frac{1}{2}) r_{1}^{2} - \sum_{i=2}^{k}(m_{1} + \frac{1}{2})r_{1}r_{i}$,
then

$\mathcal{A} = \sum_{i=2}^{k}(m_{i} - m_{1} -\frac{1}{2}) r_{1} r_{i}$ $\hspace{0.4cm}$ where $m_{i} > m_{1} + 1$
$$
$$
Note that, with our hypothesis that $E\mid_{P}$ is not pure, we have $m_i\geq m_1+1$ so
$A > 0$.

\begin{proposition}
\mylabel{propdeltalocprime}
If $E|_P$ is not pure, then let $E'$ be the elementary transformation considered above. The local invariant
satisfies
$$
\Delta^{Vb}_{loc}(E', P) = \Delta^{Vb}_{loc}(E,P) - \frac{1}{r} \sum_{i=2}^{k}(m_{i} - m_{1} -\frac{1}{2}) r_{1} r_{i}
 \hspace{0.3cm} where \hspace{0.3cm} m_{i} > m_{1} + 1.
 $$
In particular, $\hspace{0.3cm}$ $\Delta^{Vb}_{loc}(E', P) < \Delta^{Vb}_{loc}(E,P)$.
\end{proposition}

If $E'$ is pure then $\Delta ^{Vb}_{loc}(E',P)=0$,
if not we can continue by applying the elementary transformation process to $E'$ and so on, until the result is pure.
The resulting theorem can be viewed as a local analogue of the Bogomolov-Gieseker inequality.

\begin{theorem}
If $E$ is a vector bundle on $X$ and $P\cong \mathbf{P}^1\subset X$ is the exceptional divisor of blowing up a smooth
point $\check{P}\in \check{X}$, then $\Delta ^{Vb}_{loc}(E,P)\geq 0$, and $\Delta ^{Vb}_{loc}(E,P)=0$ if and only if
$E\cong \mu ^{\ast}(\check{E})$ is the pullback of a bundle from $\check{X}$.
\end{theorem}

The invariant $\Delta ^{Vb}_{loc}(X,P)$ also provides a bound for $m_i-m_1$.
\begin{corollary}
If $E\mid_{P} = \bigoplus_{i=1}^{k}\mathcal{O}(m_{i})^{r_{i}}$,  where $m_{1} < m_{2} < ... <m_{k}$. Then

$$
\Delta^{Vb}_{loc}(E', P) \geq 0 ;  \hspace{0.5cm} \Delta^{Vb}_{loc}(E, P) \geq \frac{1}{r} \sum_{i=2}^{k}(m_{i} - m_{1} - \frac{1}{2})r_1 r_{i}.
$$
\end{corollary}
\section{Modification of filtrations due to elementary transformations}

Given two bundles $E$ and $F$ such that $E|_U\cong F|_U$, then $F$ may be
obtained from $E$ by a sequence of elementary transformations. We therefore analyse what
happens to the filtrations along the divisor components $D_i$ different from exceptional
divisors $P_u$, in the case of an elementary transformation.

Suppose $E'$ is obtained from $E$ by an elementary transformation.

We have bundles $Gr (D_i, k; E)$ and $Gr (D_i, k; E')$ over $D_i$.
In order to follow the modification of the formula for $\Delta$ we need to
consider this change.

For the bundle $E$ we have a filtration by full flags $F^i_k\subset E|_{D_i}$.
Suppose $E|_{D_0}\rightarrow Q$ is a quotient (locally free on $D_0$)
and let $E'$ be the elementary
transformation fitting into the exact sequence
$$
0\rightarrow E' \rightarrow E \rightarrow Q \rightarrow 0.
$$

\begin{lemma}
Suppose $i\neq 0$ so $D_i\cap D_0$ is transverse.
Tensoring this exact sequence with $\mathcal{O} _{D_i}$ yields an exact sequence
$$
0\rightarrow E'|_{D_i} \rightarrow E|_{D_i} \rightarrow Q |_{D_i}\rightarrow 0.
$$
\end{lemma}
\begin{proof}
In fact we get a long exact sequence
$$
Tor^1_{\mathcal{O} _X}(\mathcal{O} _{D_i},Q)\rightarrow
E'\otimes \mathcal{O} _{D_i} \rightarrow E \otimes \mathcal{O} _{D_i}\rightarrow Q \otimes \mathcal{O} _{D_i}\rightarrow 0,
$$
but the facts that $Q$ is locally free on $D_0$ and $D_i$ is transverse to $D_0$ imply
that $Tor^1_{\mathcal{O} _X}(\mathcal{O} _{D_i},Q)=0$.
\end{proof}

This lemma says that $E'|_{D_i}$ is an elementary transformation of $E|_{D_i}$.
Notice that since $D_i$ is a curve, $Q|_{D_i}$ is a skyscraper sheaf.

Define $F'^i_k:= F^i_k\cap (E'|_{D_i})$. It is a subsheaf of $(E'|_{D_i})$.
Furthermore it is saturated, that is to say the quotient is torsion-free.
To show this,
suppose $s$ is a section of $(E'|_{D_i})$ which is contained in $F'^i_k$
over an open set. Then it may be seen as a section of $E|_{D_i}$ which
is contained in $F^i_k$ over an open set, but $F^i_k$ is saturated
so the section is contained in $F^i_k$. Hence by definition the section is
contained in $F'^i_k$.

Thus, we have defined a filtration $F'^i_k$ by sub-vector bundles. The
same argument says that $F'^i_k$ is saturated in $F'^i_{k+1}$, so the
quotients $Gr(D_i,k;E')= F'^i_{k+1}/F'^i_k$ are locally free; since they are
line bundles over the open set, they are line bundles on $D_i$.

Consider the morphism, induced by $(E'|_{D_i})\rightarrow (E|_{D_i})$:
$$
F'^i_{k+1}/F'^i_k \rightarrow F^i_{k+1}/F^i_k.
$$
By the definition of $F'^i_k$ it is seen that this morphism is an injection
of sheaves. Consider the cokernel. If $s$ is a section of $F^i_{k+1}/F^i_k$
and if $z$ is a local coordinate on $D_i$ such that $z=0$ defines the
intersection point $D_0\cap D_i$, then we claim that $zs$ must be in the
image of $F'^i_{k+1}/F'^i_k$. Lift $s$ to a section also denoted $s$ of $F^i_{k+1}$.
Then, thought of as a section of $E|_{D_i}$, notice that $zs$ projects to $0$ in
$Q|_{D_i}$. This may be seen by further extending to a section of $E$ and
extending $z$ to a coordinate function defining $D_0$; noting that $Q$ is supported
scheme-theoretically on $D_0$ so $zs$ projects to $0$ in $Q$.

From the exact sequence, we conclude that $zs$ is in the
image of $F'^i_{k+1}/F'^i_k$.  Hence, there are two cases:
\newline
(1) the map $F'^i_{k+1}/F'^i_k \rightarrow F^i_{k+1}/F^i_k$ is an isomorphism; or
\newline
(2) we have $F'^i_{k+1}/F'^i_k = F^i_{k+1}/F^i_k \otimes _{\mathcal{O} _{D_i}} \mathcal{O} _{D_i}(-D_0\cap D_i)$.

In the first case (1),
$$
c_1^{D_i}(Gr (D_i,k; E'))= c_1^{D_i}(Gr (D_i,k; E)).
$$
In the second case (2),
$$
c_1^{D_i}(Gr (D_i,k; E'))= c_1^{D_i}(Gr (D_i,k; E)) - [D_0\cap D_i].
$$
Applying $(\xi _i)_{\ast}$ gives the following proposition.

\begin{proposition}
Suppose $E'$ is an elementary transformation of $E$. Then  there exist a unique invariant $\emph{deg}_{loc}$ that satisfy the following properties:
$$
\emph{deg}_{loc}\left(D_j, k; \check{E}\right) : = 0,
$$
$$
\emph{deg}_{loc}\left(D_j, k; \check{E}(m.P_i)\right) : = m,
$$
and for divisor components $D_i$ intersecting $D_0$ transversally, the change in Chern class of the
associated-graded pieces is
$$
\emph{deg}_{loc}\left(Gr(D_j, k; E'), P_i \right) : = \emph{deg}_{loc}\left(Gr(D_j, k; E), P_i\right) - \tau (E,E'; k)
$$
where $\tau (E,E'; k) = 0$ or $1$ in cases (1) or (2) respectively.
\end{proposition}

\begin{definition}
\label{deglocdef}
Let $S$ and $\check{S}$ are two sheaves of finite length  with support at points $P_i$. Suppose We have bundles $Gr (D_j, k, E)$ and $Gr (D_j, k, \check{E})$
over $D_j$ respectively $\check{D}_j$. Let $F$ be the intersection of subsheaves $Gr (D_j, k, E)$ and $Gr (D_j, k, \check{E})$.
Define $lg$ to be the length, and let $lg(S,P_i)$ be the length of the part supported set-theoretically at $P_i$. Thus
$$
lg(S) = \sum _{i}lg (S, P_i)
$$
and similarly for $\check{S}$.
Consider the sequences:
$$
F \longrightarrow Gr(D_j, k, E) \longrightarrow S \longrightarrow 0
$$
$$
F \longrightarrow Gr(D_j, k, \check{E}) \longrightarrow \check{S} \longrightarrow 0.
$$
Define
$$
\emph{deg}_{loc} \left(Gr(D_j, k;E, P_i)\right) := lg(S,P_i) - lg(\check{S},P_i).
$$
\end{definition}
Then
$$
\emph{deg} \left(Gr(D_j, k, E)\right) = \emph{deg}(F) + lg(S)
$$
$$
\emph{deg} \left(Gr(D_j, k, \check{E})\right) = \emph{deg}(F) + lg(\check{S}),
$$
therefore
$$
lg(S) - lg(\check{S}) = \sum_{P_i}\left[lg(S, P_i) - lg(\check{S}, P_i)\right] = \sum_{P_i} \emph{deg}_{loc} \left(Gr(D_j, k;E, P_i)\right).
$$
Suppose $i \neq 0$, if $D_j$ are  non-exceptional divisor. Then
$$
\emph{deg}\left(Gr(D_j, k; E\right) = \emph{deg}\left(Gr(\check{D}_j, k; \check{E} \right) + \sum_{P_i} \emph{deg}_{loc}\left(Gr(D_j, k, E), P_i)\right).
$$
This completes the proof of the proposition.
\section{The local parabolic invariant}

Let $E$ be a bundle, with $\beta(D_0, k) = 0,  \forall k$. Then we would like
to define the terms in the following equation:
$$
\Delta^{Par}(E) = \Delta^{Par}(\check{E})  +\sum_{P_{i}}\Delta^{Par}_{loc}(E, P_i).
$$

Assume that $\check{D}$ is a union of smooth divisors meeting in some multiple points.
The divisor $D$ is obtained by blowing up the points $\check{P}_u$ of multiplicity $\geq 3$.
Let
$$
\varphi : X\rightarrow \check{X}
$$
be the birational transformation.
We use the previous formula to break down $\Delta ^{Par}(E)$ into a global contribution
which depends only on $\check{E}$, plus a sum of local contributions depending
on the choice of extension of the parabolic structure across $P_u$.

Let $\check{\mathcal{S}}$ denote the set of divisor components in $\check{D}$ (before blowing-up)
and
define the global term $\Delta ^{Par}(\check{E})$ by the formula
$$
$$
$\Delta^{Par}(\check{E}) := \Delta^{Vb}(\check{E}) $
$$
$$
$\hspace{2.2cm} + \  \sum_{i \in \check{\mathcal{S}}}\sum_{k=1}^r \beta (\check{D}_i,k).\emph{deg}\left(Gr(\check{D}_i, k)\right) $
$$
$$
$\hspace{2.2cm} - \ \frac{1}{2} \ \sum_{i\in \mathcal{S}}\sum_{k=1}^r
\beta (\check{D}_i,k)^{2}.[\check{D}_{i}]^{2}$
$$
$$
$\hspace{2.2cm} - \ \frac{1}{2}\sum_{i \neq j}
\sum_{k=1}^r  \sum_{y \in Irr(\check{D}_{i} \cap \check{D}_{j})} \beta (\check{D}_i,k).\beta (\check{D}_j,\sigma (y,i,j)(k)).[y].$
$$
$$
This formula imitates the formula for $\Delta ^{Par}$ by considering only pairwise intersections
of divisor components even though several different pairwise intersections could occur at the
same point.
Recall that $[D_i]^2 = [\check{D}_i]^2 - m$ where $m$ is the number of points on $\check{D}_i$
which are blown up to pass to $D_i$.

To define the local terms, suppose at least one of the divisors, say $D_0=P$, is the exceptional locus
for a birational transformation blowing up the point $\check{P}$.
We define a local contribution $\Delta ^{Par}_{loc}(E,P)$ to $\Delta ^{Par}$ by isolating
the local contributions in the previous formula.

Notice first of all that for any $D_i$ meeting $P$ transversally,
we have defined above $\emph{deg}_{loc}\left(Gr(E;D_i, k), P\right)$,
the local contribution at $P$,
in such a way that
$$
\emph{deg}\left(Gr(E;D_i, k)\right)
= \emph{deg}\left(Gr(\check{E}; D_i, k)\right) + \sum _{P_u} \emph{deg}_{loc}\left(Gr(E;D_i, k), P_u\right)
$$
where the sum is over the exceptional divisors $P_u$ meeting $D_i$, which correspond to the
points $\check{P}_u\in \check{D}_i$ which are blown up.

Let $\mathcal{S}(P)$ denote the set of divisor components which meet $P$ but not including
$P=D_0$ itself. Define
$$
$$
$\Delta^{Par}_{loc}(E,P) := \Delta^{Vb}_{loc}(E,P) $
$$
$$
$\hspace{2.2cm} + \  \sum_{k=1}^r \beta (P,k).
\emph{deg}\left(Gr(E;P, k)\right) $
$$
$$
$\hspace{2.2cm} + \  \sum_{i \in \mathcal{S}(P)}\sum_{k=1}^r \beta (D_i,k).
\emph{deg}_{loc}\left(Gr(E;D_i, k), P\right) $
$$
$$
$\hspace{2.2cm} + \ \frac{1}{2} \ \sum_{k=1}^r
\beta (P,k)^{2}$
$$
$$
$\hspace{2.2cm} + \ \frac{1}{2} \ \sum_{i\in \mathcal{S}(P)}\sum_{k=1}^r
\beta (D_i,k)^{2}$
$$
$$
$\hspace{2.2cm} - \ \sum_{i \in  \mathcal{S}(P)}
\sum_{k=1}^r  \beta (D_i,k).\beta (P,\sigma (i,P)(k)).[y]$
$$
$$
$\hspace{2.2cm} +\ \frac{1}{2}\sum_{i\neq j,  \check{P}\in \check{D}_i\cap \check{D}_j}
\sum_{k=1}^r  \beta (\check{D}_i,k).\beta (\check{D}_j,\sigma (\check{P},i,j)(k)).[\check{P}]$
$$
$$

In the next to last term,
$\sigma (i,P):=\sigma (y,i,v)$ where $P=D_v$ and $y$ is the unique intersection
point of $P=D_v$ and $D_i$. The factor of $1/2$ disappears because we are
implicitly choosing an ordering of the indices $i,j=0$ which occur here.
The last term is put in to cancel with the
corresponding term in the global expression for
$\check{E}$ above, and $[\check{P}]$ designates any lifting of the point
$\check{P}$ to a point on $P$.

\begin{theorem}
With the above definitions, we have
$$
\Delta ^{Par}(E) = \Delta ^{Par}(\check{E}) + \sum _{P_u} \Delta ^{Par}_{loc}(E, P_u),
$$
where the sum is over the exceptional divisors.
\end{theorem}
\begin{proof}
This follows by comparing the above definitions with the formula of Proposition \ref{propepar}.
\end{proof}

Let $\varphi ^{\ast}\check{E}$ denote the parabolic bundle on $X$ given by using
the trivial extension $\varphi^{\ast}E$ as underlying vector bundle, and setting
$\beta (P_u,k):= 0$ for all exceptional divisor components $P_u$.
Note that $\Delta^{Vb}_{loc}(\varphi^{\ast}E,P)=0$.
Then
\\
$\Delta^{Par}_{loc}(\varphi ^{\ast}\check{E},P) =  \ \frac{1}{2} \ \sum_{i\in \mathcal{S}(P)}\sum_{k=1}^r
\beta (\check{D}_i,k)^{2}$
$$
$$
$\hspace{2.2cm} +\ \frac{1}{2}\sum_{i\neq j,  \check{P}\in \check{D}_i\cap \check{D}_j}
\sum_{k=1}^r  \beta (\check{D}_i,k).\beta (\check{D}_j,\sigma (\check{P},i,j)(k)).[\check{P}]$,
\\
and
$$
$$
$\Delta ^{Par}_{loc}(E,P)-
\Delta ^{Par}_{loc}(\varphi ^{\ast}\check{E},P) = \Delta^{Vb}_{loc}(E,P)$
$$
$$
$\hspace{2.2cm} + \  \sum_{k=1}^r \beta (P,k).
\emph{deg}\left(Gr(E;P, k)\right) $
$$
$$
$\hspace{2.2cm} + \  \sum_{i \in \mathcal{S}(P)}\sum_{k=1}^r \beta (D_i,k).
\emph{deg}_{loc}\left(Gr(E;D_i, k), P\right) $
$$
$$
$\hspace{2.2cm} + \ \frac{1}{2} \ \sum_{k=1}^r
\beta (P,k)^{2}$
$$
$$
$\hspace{2.2cm} - \ \sum_{i \in  \mathcal{S}(P)}
\sum_{k=1}^r  \beta (D_i,k).\beta (P,\sigma (i,P)(k)).[y]$.

A different local-global decomposition may be obtained by noting that
$$
\Delta ^{Par}(E)=\Delta^{Par}(\varphi ^{\ast}\check{E}) +
\sum _u \left(
\Delta ^{Par}_{loc}(E,P_u)-
\Delta ^{Par}_{loc}(\varphi ^{\ast}\check{E},P_u)
\right)
$$
with the local terms $(\Delta ^{Par}_{loc}(E,P_u)-
\Delta ^{Par}_{loc}(\varphi ^{\ast}\check{E},P_u))$ being given by the previous formula.
\section{Normalization via standard elementary transformations}

There is another modification of parabolic structures due to elementary transformations.
This may also be viewed as a shift of the parabolic structures in the viewpoint of
a collection of sheaves. If $E_{\cdot}=\{ E_{\alpha _1,\ldots , \alpha _n}\}$
is a parabolic sheaf, then we can shift the filtration at the $i$-th place
defined by
$$
(C^i_{\theta}E)_{\alpha _1,\ldots , \alpha _n}:= E_{\alpha _1,\ldots , \alpha _i-\theta , \ldots , \alpha _n}.
$$
This may also be viewed as tensoring with a parabolic line bundle
$$
C^iE = E\otimes \mathcal{O} (\theta D_i).
$$
The weights of the parabolic structure $C^i_{\theta}E$ along $D_i$
are of the form $\alpha _i+\theta$ for $\alpha _i$ weights of $E$.

In the point of view of a vector bundle with filtration, it may correspond to doing an
elementary transformation. Suppose $0 <\theta < 1$. Then
$$
(C^i_{\theta}E)_0 = E_{0,\ldots , 0,-\theta , 0,\ldots ,0}
$$
and we have an exact sequence
$$
0\rightarrow (C^i_{\theta}E)_0 \rightarrow E_0 \rightarrow (E_0 / F^i_{-\theta}E_0)\rightarrow 0.
$$
Therefore $(C^i_{\theta}E)_0$ is obtained by elementary transformation of $E_0$ along one
of the elements of the parabolic filtration on the divisor $D_i$.

This is specially useful when the rank is $2$.
Suppose $rk(E)=2$. There is a single choice for the elementary transformation.
If the weights of $E$ at $D_i$ are $\alpha^{tot} _i-\beta _i$ and $\alpha^{tot}_i + \beta _i$
then the weights of the elementary transformation will be $\theta  +\alpha^{tot} _i+\beta _i-1$ and $\theta  +\alpha^{tot} _i-\beta _i$.
The shift $\theta $ should be chosen so that these lie in $(-1,0]$.
The new weights may be written as
$$
(\tilde{\alpha}_i^{tot} -\tilde{\beta}_i,\quad \tilde{\alpha}_i^{tot} + \tilde{\beta}_i
$$
with
$$
\tilde{\alpha}_i^{tot}:= (\theta  +\alpha^{tot} _i-\frac{1}{2})
$$
which is the new average value, and
$$
\tilde{\beta}_i:=\frac{1}{2} -\beta _i.
$$
\begin{corollary}
\label{quarter}
In the case $rk(E)=2$,
by replacing $E$ with its shift $C^i_{\theta}E$ if necessary, we may assume that
$$
0\leq \beta _i \leq \frac{1}{4}.
$$
\end{corollary}
\begin{proof}
If $\beta _i> \frac{1}{4}$ then do the shift which corresponds to an elementary
transformation; for the new parabolic structure $\tilde{\beta}_i=\frac{1}{2}-\beta _i$
and $0\leq \tilde{\beta}_i \leq \frac{1}{4}$.
\end{proof}
\section{The rank two case}

In order to simplify the further constructions and computations, we now restrict
to the case when $E$ has rank 2.
The parabolic structures along $D_i$ are rank one subbundles $F^i\subset E|_{D_i}$.
The associated graded pieces are $Gr (D_i,1)=F^i$ and $Gr(D_i, 2)=E|_{D_i}/F^i$.
The normalized weights may be written as
$$
\beta (D_i,1)=-\beta _i, \quad \beta (D_i, 2) =\beta _i
$$
with $0\leq \beta _i<\frac{1}{2}$, and by Corollary \ref{quarter} we may furthermore suppose $0\leq \beta _i\leq\frac{1}{4}$.

Define $\deg ^{\delta}(E_{D_i},F^i) := \left( \deg (E_{D_i}/F^i)-\deg (F^i)\right)$.
This has a local version as discussed in Definition \ref{deglocdef},
$$
\deg^{\delta}_{loc}(E_{D_i},F^i,P) := \left( \deg _{loc}(E_{D_i}/F^i,P)-\deg _{loc}(F^i,P)\right)
$$
whenever $D_i$ meets $P$ transversally.

The main formula may now be rewritten:

\begin{proposition}
$\Delta^{Par}(E) = \Delta^{Vb}(E) $
$$
$$
$\hspace{2.2cm} + \  \sum_{i \in \mathcal{S}} \beta _i\deg^{\delta}(E_{D_i},F^i)  $
$$
$$
$\hspace{2.2cm} - \  \ \sum_{i\in \mathcal{S}}
\beta _i^{2}.[D_{i}]^{2}$
$$
$$
$\hspace{2.2cm} - \ \sum_{i \neq j}
\sum_{y \in Irr(D_{i} \cap D_{j})} \tau (y,i,j)\beta _i\beta _j.[y]$
\\
where $\tau (y,i,j)=1$ if $F^i(y)=F^j(y)$ and $\tau (y,i,j)=-1$ if $F^i(y)\neq F^j(y)$
as subspaces of $E(y)$.

Similarly for the local parabolic invariants, denoting $P=D_0$ we have
$$
$$
$\Delta ^{Par}_{loc}(E,P)-
\Delta ^{Par}_{loc}(\varphi ^{\ast}\check{E},P) = \Delta^{Vb}_{loc}(E,P)$
$$
$$
$\hspace{2.2cm} + \  \beta _0.
\deg ^{\delta}(E_{D_0},F^0) $
$$
$$
$\hspace{2.2cm} + \  \sum_{i \in \mathcal{S}(P)} \beta _i.
\deg ^{\delta}_{loc}(E_{D_i},F^i,P)  $
$$
$$
$\hspace{2.2cm} + \
\beta _{0}^{2}$
$$
$$
$\hspace{2.2cm} - \ 2 \sum_{i \in  \mathcal{S}(P)}
\tau (i,P) \beta _i\beta _0 .[y]$.
\end{proposition}
\begin{example}
Let $E$ be a non pure rank two bundle, we have $E\mid_{P} = \mathcal{O}(m_{1}) \oplus \mathcal{O}(m_{2})$, then
$\hspace{0.4cm}$ $\Delta^{Vb}_{loc}(E,P) \geq \frac{1}{2}(m_{2} - m_{1} - \frac{1}{2})$.
\end{example}
\begin{example}
Let $E$ be a rank 2 bundle with $E|_{P} = \mathcal{O} \oplus \mathcal{O}(1)$. The reduction by elementary transformation
is $E|_{P} = \mathcal{O} \oplus \mathcal{O}(1) \rightsquigarrow E'$.
We get an exact sequence
$$
0 \longrightarrow \mathcal{O}(1) \longrightarrow E'|_{P} \longrightarrow \mathcal{O}(1) \longrightarrow 0
$$
then $E'$ is pure, $E' = \mu^{\star}\check{E}(-P)$.
\end{example}

Suppose we start with the bundle $E$, by doing the sequence of elementary transformation  we get $\check{E}(m.P_i)$. Number the sequence in opposite direction,
we get a sequence of bundles of the form:
$$
\check{E}(m.P_i) = E(0),\hspace{0.1cm} E(1),\hspace{0.1cm} E(2),\hspace{0.1cm} ...,\hspace{0.1cm} E(g) = E
$$
where $g$ is the number of steps, and  $E(j - 1) = (E(j))'$ for $j = 1, ...,g$.
we recall that if $E|_P \cong \mathcal{O}(m_1) \oplus \mathcal{O}(m_2)$ with $m_1  \leq m_2$ then $\mu(E) = m_2 - m_ 1$.
Also  $\mu(E) = 0 \Longrightarrow E = \check{E}(m_iP_i)$.
Furthermore
if $m_1 < m_2$ then $\mu(E') < \mu(E)$.
We see that
$$
0 = \mu(E_0) < \mu(E_1) < \mu(E_2) < ... < \mu(E_{g - 1}).
$$

To calculate $\Delta(E, P)_{loc}$  we use the proposition \ref{propdeltalocprime} applied to each $E(j)$:
$$
\Delta^{Vb}_{loc}(E(j)', P_0) = \Delta^{Vb}_{loc}(E(j), P_0) - \frac{1}{2} \mathcal{A},
$$
where
$$
\mathcal{A} = \sum_{i = 2}^{k}(m_2 - m_1 - \frac{1}{2})r_1r_i = m_2(E(j)) - m_1(E(j)) - \frac{1}{2} = \mu(E(j)) - \frac{1}{2}.
$$
Therefore
$$
\Delta^{Vb}_{loc}(E(j - 1)) = \Delta^{Vb}_{loc}(E(j)) - \frac{1}{2}(\mu(E(j)) - \frac{1}{2}),
$$
and putting them all together,
$$
\Delta^{Vb}_{loc}(E(0)) = 0; \hspace{0.2cm} \Delta^{Vb}_{loc}(E(g)) = \frac{1}{2}\sum_{i = 1}^{g}(\mu(E(j)) - \frac{1}{2}).
$$
We have
\begin{equation}
\mylabel{varphi}
\mu(E(j-1)) < \mu((E(j)),
\end{equation}
so  $\mu(E(j)) \geq j$.

Now we divide the work in two parts,  first term $\mu(E(g)) - \frac{1}{2}$, then the sum of the others.
For each $\mu(E(j))$ is at least one  greater than the previous one, this gives that the sum of the other terms is  at least equal to
$ (1 + 2 + 3 + ... + (g-1)) - \frac{(g-1)}{2}$.
Then
$$
\Delta^{Vb}_{loc}(E(g)) \geq \frac{1}{2}(1 + 2 + 3 + ... + (g-1)) - \frac{g}{4} + \frac{1}{2}\mu(E(g))
$$
$$
 = \frac{g(g-1)}{4} - \frac{g}{4} + \frac{1}{2}\mu(E).
$$
We have therefore proven the following:
\begin{proposition}
\mylabel{deltavbloc}
If $E$ is a bundle which is brought to pure form in $g \geq 1$ steps of elementary transformation, and $\mu(E) = m_2(E|_P) - m_1(E|_P)$, then
we have the lower bound
$$
\Delta^{Vb}_{loc}(E) = \Delta^{Vb}_{loc}(E(g)) \geq \frac{g^2 - 2g}{4} + \frac{1}{2}\mu (E).
$$
\end{proposition}

We have for each $1\leq k \leq g$,
$$
\left| \deg ^{\delta}_{loc}(E(k), F^i,P) - \deg ^{\delta}_{loc}(E(k-1), F^i,P) \right| \leq 1,
$$
but also $E(0)=\varphi ^{\ast}(\check{E})$ and $\deg ^{\delta}_{loc}(E(0), F^i,P) =0$, so
$$
\left| \deg ^{\delta}_{loc}(E(g), F^i,P)  \right| \leq g.
$$

Also along $P$ we have $E_P=\mathcal{O}(m_1)\oplus \mathcal{O}(m_2)$ with $m_1\leq m_2$. For any subbundle $F^0\subset E_P$
we have $\deg (F^0) \leq m_2$ and $\deg (E_P/F^0)\geq m_1$ so
$$
\deg ^{\delta}(E_P, F^0)\geq m_1-m_2 = -\mu (E).
$$
Then
$$
$$
$\Delta ^{Par}_{loc}(E,P)-
\Delta ^{Par}_{loc}(\varphi ^{\ast}\check{E},P) = \Delta^{Vb}_{loc}(E,P)$
$$
$$
$\hspace{2.2cm} + \  \beta _0.
\deg ^{\delta}(E_{D_0},F^0) $
$$
$$
$\hspace{2.2cm} + \  \sum_{i \in \mathcal{S}(P)} \beta _i.
\deg ^{\delta}_{loc}(E_{D_i},F^i,P)  $
$$
$$
$\hspace{2.2cm} + \
\beta_{0}^{2}$
$$
$$
$\hspace{2.2cm} - \ 2 \sum_{i \in  \mathcal{S}(P)}
\tau (i,P) \beta _i\beta _0 .[y]$
$$
$$
$\hspace{2.2cm} \geq\frac{g^{2} - 2g}{4} \hspace{0.2cm} + \hspace{0.2cm} \frac{1}{2}\mu (E)$
$$
$$
$\hspace{2.2cm} - \  \beta _0\mu (E)$
$$
$$
$\hspace{2.2cm} - \sum _{i\in  \mathcal{S}(P)}\beta _i g$
$$
$$
$\hspace{2.2cm} + \
\beta _{0}^{2}$
$$
$$
$\hspace{2.2cm} - \ 2 \sum_{i \in  \mathcal{S}(P)}
\beta _i\beta _0 .$

But we know that $\mid\beta_i\mid \leq \frac{1}{2}$. Then we get the following theorem.
\begin{theorem}
$$
\Delta^{Par}_{loc}(E; P) \geq \Delta ^{Par}_{loc}(\varphi ^{\ast}\check{E},P) + \frac{g^{2} - 2g}{4} - \frac{g +1}{2}.\kappa  .
$$
Where $\kappa = \# \mathcal{S}(P)$  is  the number of divisors of $D_i$ meeting $P$.
\end{theorem}
\begin{theorem}
If $\check{E}$ is a vector bundle of rank $2$ on $\check{X}$ with parabolic structures
on the components $\check{D}_i$, then on $X$ obtained by blowing up the multiple points of
$\check{D}$, the parabolic invariant $\Delta ^{Par}_{loc}(E, P)$ attains a minimum for
some extension of the bundle $E$ and some parabolic structures on the exceptional loci.
\end{theorem}
\begin{proof}
From the above theorem, the number of elementary transformations $g$ needed to get to any
$E$ with $\Delta^{Par}_{loc}(E; P) \leq \Delta ^{Par}_{loc}(\varphi ^{\ast}\check{E},P)$, is bounded.
Furthermore the number of numerical possibilities for the degrees
$\deg ^{\delta}_{loc}(E_{D_i},F^i,P) $ and $\deg ^{\delta}(E_{P},F^0,P) $ leading to such a minimum, is finite.
The parabolic weight $\beta _0$ may be chosen to lie in the closed interval $[0,\frac{1}{4}]$,
so the set of possible numerical values lies in a compact subset; hence a minimum is attained.
\end{proof}

Denote the parabolic extension which achieves the minimum by $E^{min}$.
There might be several possibilities, although we conjecture that usually it is unique.
Thus
$$
\Delta ^{Par}_{loc}(E^{min}, P) = {\rm min}_E \left(  \Delta ^{Par}_{loc}(E, P)\right) .
$$
With the minimum taken over all parabolic extensions $E$ of $\check{E}|_{\check{U}}$ across
the exceptional divisor $P$.

The minimal $E^{min}$
exists at each exceptional divisor and they fit together to give a global parabolic bundle.
Define
$$
$$
$\Delta ^{Par}_{min}(\check{E}):= \Delta ^{Par}(E^{min})$
$$
$$
$\hspace{2.2cm}
= \Delta ^{Par}(\check{E}) + \sum _{P_u} \Delta ^{Par}_{loc}(E^{min},P_u).$
\subsection{Panov differentiation}

D. Panov in his thesis \cite{Panov} used the idea of differentiation with respect
to the parabolic weight. A version of this
technique allows us to gain more precise information on the minimum.

\begin{lemma}
Let $E=E^{min}$ be the parabolic bundle extending $\check{E}|_{\check{U}}$
which achieves the minimum value $\Delta ^{Par}_{loc}(E^{min},P)$.
By making an elementary transformation we may assume $0\leq \beta _0\leq \frac{1}{4}$.
Denote also by $E$
the underlying vector bundle.
Then for any subbundle $F'\subset E|_P$ we have
$$
\deg (E|_P/F')- \deg (F')   \geq -\kappa .
$$
Thus if $E|_P = \mathcal{O}_P(m_1)\oplus \mathcal{O}_P(m_2)$ then
$$
|m_2-m_1| \leq \kappa .
$$
Where $\kappa = \# \mathcal{S}(P)$  is  the number of divisors of $D_i$ meeting $P$.
\end{lemma}
\begin{proof}
We show that $\deg (E|_P/F')- \deg (F')   \geq -\kappa -\frac{1}{4}$ which implies the
stated inequality since the left side and $\kappa$ are integers.
Let $F=F^0\subset E|_P$ be the subbundle corresponding to the parabolic structure $E^{min}$.
Consider two cases:
\\
(i) if $F^0$ is the destabilizing bundle of $E|_P$ and $\beta _0 >0$; or
\\
(ii) if $F^0$ is not the destabilizing bundle of $E|_P$, or else $\beta _0=0$.

In case (i) note that $\beta_0$ may be allowed to range in the full interval $[0,\frac{1}{2})$ so
the invariant $\Delta ^{Par}_{loc}(E,P)$ is a local minimum considered as a function of $\beta_0\in (0,\frac{1}{2})$.
Then
$$
\frac{d}{d\beta _0} \Delta ^{Par}_{loc}(E,P) = 0.
$$
This gives the formula
$$
$$
$
\deg ^{\delta}(E_{D_0},F^0)
=2 \
\beta _0 + \ 2 \sum_{i \in  \mathcal{S}(P)}
\tau (F^i,F^0) \beta _i $,
so
$$
$$
$\deg ^{\delta}(E_{D_0},F^0) \geq - \frac{\kappa}{2}$.

Since $F^0$ is the destabilizing bundle it implies that
$$
$$
$\deg ^{\delta}(E_{D_0},F') \geq -\frac{\kappa}{2}$
\\
for any other subbundle $F'$ also, which is stronger than the desired inequality in this case.

In case (ii) we have $\beta _0.\deg ^{\delta}(E_{D_0},F^0) \geq 0$ because in the contrary case
that would imply that $F^0$ is the destabilizing subbundle.
Suppose $F'\subset E|_P$ is a possibly different subbundle such that
$$
\deg (E|_P/F')- \deg (F')   < -\frac{1}{4}(1+4\kappa) .
$$
Then make a new parabolic structure $E'$ using $F'$ instead of $F$, with parabolic weight $\beta '_0 =
\frac{1}{4}$. We have
$$
$$
$\Delta ^{Par}_{loc}(E',P)-
\Delta ^{Par}_{loc}(E,P) = $
$$
$$
$\hspace{2.2cm}  \   \frac{1}{4}
\deg ^{\delta}(E_{D_0},F') $
$$
$$
$\hspace{2.2cm} - \  \beta _0.
\deg ^{\delta}(E_{D_0},F^0) $
$$
$$
$\hspace{2.2cm} + \frac{1}{16} $
$$
$$
$\hspace{2.2cm} - \
\beta _0^{2}$
$$
$$
$\hspace{2.2cm} + \ 2 \sum_{i \in  \mathcal{S}(P)}
\tau (F,F^i) \beta _i\beta _0 .[y]$
$$
$$
$\hspace{2.2cm} - \ \frac{1}{2}\sum_{i \in  \mathcal{S}(P)}
\tau (F',F^i) \beta _i .[y]$
$$
$$
$\leq  \frac{1}{4}
\deg ^{\delta}(E_{D_0},F')+ \frac{1}{16}(1+4\kappa) $
$$
$$
$<0.$

This contradicts minimality of $E^{min}$, which shows the desired inequality.
\end{proof}

\begin{corollary}
In the case of $3$ divisor components $\kappa = 3$ and the minimal extension $E^{min}$
satisfies $|m_2-m_1| \leq 3$. It is connected to $\varphi ^{\ast}(\check{E})$ by at most three
elementary transformations.
\end{corollary}

This should permit an explicit description of all possible cases for $\kappa = 3$, we start on this below.
\subsection{The Bogomolov-Gieseker inequality}

Suppose $C\subset \check{X}$ is an ample curve meeting $\check{D}$ transversally. Then
$\check{E}|_C$ is a parabolic bundle on $C$.

\begin{proposition}
Suppose $\check{E}|_C$ is a stable parabolic bundle. Then for any extension $E$ to a
parabolic bundle over $X$, there exists an ample divisor $H$ on $X$ such that
$E$ is $H$-stable. Hence $\Delta ^{Par}(E)\geq 0$. In particular $\Delta ^{Par}(E^{min}) \geq 0$.
If $\check{E}$ comes from an irreducible unitary representation of $\pi _1(\check{X}-\check{D})$
then the parabolic extension on $X$ corresponding to the same unitary representation
must be some choice of $E^{min}$.
\end{proposition}
\begin{proof}
Fix an ample divisor $H'$. Then any divisor of the form $H=nC+H'$ is ample on $X$,
and for $n$ sufficiently large $E$ will be $H$-stable. The Bogomolov-Gieseker inequality
for parabolic bundles says that $\Delta ^{Par}(E)\geq 0$ with equality if and only if
$E$ comes from a unitary representation. However, $\Delta ^{Par}(E)\geq \Delta ^{Par}(E^{min})\geq 0$
and if $E$ comes from a unitary representation then
$\Delta ^{Par}(E)= \Delta ^{Par}(E^{min})= 0$. It follows in this case that $E$ is one of the
choices of $E^{min}$.
\end{proof}



Laboratoire J. A. Dieudonn\'e, UMR 6621 Universit\'e  Nice Sophia-Antipolis,\\
06108 Nice, Cedex 2, France.\\
E-mail address: \textbf{taher@math.unice.fr, chaditaher@hotmail.fr}\\


\end{document}